\documentclass[letterpaper, 10 pt, conference]{ieeeconf}  

\IEEEoverridecommandlockouts                              

\newcommand{\hide}[1]{}

\usepackage{hyperref}

\usepackage{cite}
\usepackage{amsmath,amssymb,amsfonts,mathtools}

\usepackage{algorithmic}
\usepackage{graphicx}
\usepackage{textcomp}
\usepackage{xcolor}
\usepackage{float}
\usepackage{subcaption}
\usepackage{nicefrac}
\usepackage{mathrsfs}

\usepackage{bigints}
\newcommand{\eop}{{\hfill $\Box$}}
\usepackage{dsfont}
\usepackage{bm}
\newtheorem{thm}{Theorem}
\usepackage{hyperref}
\newtheorem{lem}{Lemma}

\newcommand{\TR}{\ifthenelse{1>2}}

\newcommand{\spio}{\eta_1 ({\bm \varrho})}
\newcommand{\spit}{\eta_2 ({\bm \varrho})}

\def\BibTeX{{\rm B\kern-.05em{\sc i\kern-.025em b}\kern-.08em
    T\kern-.1667em\lower.7ex\hbox{E}\kern-.125emX}}

\newcommand{\Walk}{S}
\newcommand{\walk}{s}
\newcommand{\IA}{A}

\newcommand{\policy}{{\pi}}
\newcommand{\policyex}{{\pi_{ex}}}

\title{\Large \bf
Stability of Polling Systems for a Large Class of   Markovian Switching Policies\thanks{The authors are listed in alphabetical order.
This work is financially supported by DST-Inria CEFIPRA LION Project. The paper is accepted at IEEE CDC 2025.}
}

\author{Konstantin Avrachenkov$^{1}$, Kousik Das$^{2}$, Veeraruna Kavitha$^{2}$ and Vartika Singh$^{3}$
\thanks{$^{1}$Inria Centre at Universit\'e C\^ote d’Azur, Sophia Antipolis, France
        {\tt\small k.avrachenkov@inria.fr}}%
\thanks{$^{2}$Department of IEOR, IIT Bombay, Mumbai, India 
        {\tt\small \{20004007, vkavitha\}@iitb.ac.in}}
\thanks{$^{3}$Vartika Singh was with University of Colorado, Colorado Springs, United States when this research was conducted  
{\tt\small singh.vsvartika@gmail.com}}
        }

\begin{document}

\maketitle

\begin{abstract}%
We consider a polling system with two queues, where  a single server is attending the queues in a cyclic order  and requires  non-zero switching times to switch between the queues.
Our aim  is to identify a fairly general and comprehensive class of Markovian
switching policies that renders the system stable. Potentially a class of policies that can cover the
Pareto frontier related to individual-queue-centric performance measures like  the stationary expected number of waiting customers in each queue;  for instance, such a  class of policies is identified recently for a  polling system near the fluid regime (with large arrival and departure rates), and we aim to include that class.  We also aim to include a second class that  facilitates switching between the queues at the instance the occupancy in the opposite queue crosses a threshold and when that in the visiting queue is below a threshold (this inclusion facilitates design of `robust' polling systems).
Towards this, we   consider a class of two-phase switching policies, 
which includes the above mentioned classes.   In  the  maximum generality,  our policies can be represented by eight parameters, while two  parameters are sufficient to represent the aforementioned  classes. 
We provide simple conditions to identify the sub-class of switching policies that ensure  system stability. 
By numerically tuning the parameters of the proposed class,  we illustrate that the proposed class can cover the Pareto frontier for the stationary expected number of customers in the two queues.
\end{abstract}


\section{Introduction}

Polling systems are a special class of queueing systems where a single server cyclically visits multiple queues, potentially incurring non-zero switching or walking times between the queues (see \cite{boxma2018}, \cite{czerniak2009fluid}, \cite{takagi}). These switching times introduce ``vacations," making polling systems non-work-conserving.  Polling systems are useful for modeling a variety of applications like, traffic signals, manufacturing systems, logistics, railway crossings etc.  (see \cite{pollingapp}, \cite{levy1990polling}, \cite{boxma2015}, \cite{ctspolling}, \cite{mixedpolling}). The analytical complexity of the polling systems heavily depends on the service discipline adopted (e.g., exhaustive or gated, see \cite{boxma2009},\cite{boxma2011}).
Tools like pseudo-conservation laws have proven to be invaluable in analyzing the total workload of such systems, providing insights into system-wide performance metrics (see \cite{boxma1987pseudo} and references therein). However, the performance of the individual  queues such as stationary average waiting times (or stationary average occupancy levels) remain difficult to be computed.

To obtain individual-queue-centric analysis (like achievable region of all possible stationary average occupancy levels in the individual queues), one needs to consider a large class of state dependent or Markovian switching policies; for example, Kleinrock  in \cite{kleinrock1967time} consider Markov priority based switching policies in the context of queueing systems without vacations. 
 For such work conserving systems, stability is ensured once the load factor is less than one. 
However,  with   non-zero switching times,  that is for non-work-conserving systems,  stability  can further depend upon the parameters of the switching policies. In other words,    stability can become  sensitive to the  switching policies, and thus it is important to identify the subclass (among any given class) of policies that are stable.

To the best of our knowledge, the focus of the polling system  based  literature     is on the  system performance. Further, except for a few   strands   (like \cite{WiOpt}, \cite{konstantin},  \cite{Perel2018}), the switching policies considered in the current literature are not exactly Markovian (specifically, the switching decisions do not depend on elaborate details such as queue lengths of the individual queues  at the time of decision).  For instance, in  \cite{boxma2015,boxma2011,pollingapp,boxma2018,ctspolling}, the authors consider variants of exhaustive, gated, or, globally gated disciplines at individual queues to determine  when to switch from the visiting queue; some of them (see \cite{boxma2015})  choose the  next queue to be visited   using more sophisticated rules (instead of traditional  cyclic routing). However, they do not consider Markovian switching rules  in the strict sense.

 The  authors in \cite{konstantin}  consider  Markovian (occupancy-driven) switching rules; however, they restrict to  finite buffer systems which are automatically stable, whereas we consider infinite buffer capacity systems. 
  In \cite{WiOpt}, the authors study Pareto frontier of polling systems  under a class of Markovian switching policies; however,  they focus on fluid polling systems. 
In \cite{Perel2018}, the authors  consider threshold-based switching policies (based on the occupancy of the unattended queue). However, they   
assume that one of the 
  queues have  finite buffer capacity and that the switching times are
exponentially distributed.  Due to  these assumptions,
they have a continuous time Markov jump process which
is  stable, once the combined load factor is strictly less
than one.  Essentially, the threshold parameters do not  impact   the stability of the  system due to  finite buffer considerations. We do not make any such assumptions  and consider  a big class of Markovian two-phase switching policies. In our case, the system stability further depends upon the parameters of the switching policy (apart from the load factor being strictly less than one).

In all, our work is a step towards   individual-queue-centric considerations by defining  a sufficiently large class of two-phase Markovian policies, for systems  
with infinite buffers at both the queues. Our switching policy is characterized by affine functions of queue-occupancy in both the queues 
and is parameterized  by eight coefficients.  The affine functions are piecewise constant, changing only at arrival or departure epochs, thereby allowing potential switching decisions at these precise epochs. 
The proposed class of  policies includes all the above mentioned classes   (e.g., threshold policies of \cite{konstantin, Perel2018}, priority-factor based policies of \cite{WiOpt}). We also numerically illustrate  that our policies can cover the Pareto frontier,   defined in terms of the stationary expected numbers in the individual queues.

The main aim  of this initial work is to identify the subclass of proposed   switching policies that are stable  (further theoretical investigation with respect to control and optimization is planned for future exploration).
The stochastic process representing the system evolution in our case is not even Markovian.   The idea   is to  study  the stability of the  chain,  observed at the server-visit-epochs to one of the queues. For Poisson arrivals, this forms   a discrete-time Markov chain. Using  the  Foster-Lyapunov stability criterion,  we show that a single inequality, defined in terms of two  parameters of the switching policy  determines the stability. We then prove   the time stationarity of  the system    by defining appropriate  regeneration epochs and by showing that the expected value of the corresponding regeneration time is finite; the latter statement  follows from the stability of the discrete chain.

\section{Polling System and Switching Policies}
We consider a polling system with two queues ${Q}_1$ and ${Q}_2$, where a single server attends the customers of the two queues in a cyclic order and switches between the queues with non-zero switching times. The arrival of customers to $Q_i$ (for $i\in \{1,2\}$)  is modeled as a Poisson process with rate $\lambda_i$, independent of that in the other queue indexed by $j:=3-i$ (when $i=1$, $j=2$ and when $i=2$, $j=1$).  The queues have infinite buffer capacity and follow
 the non-resume preemptive service discipline;  basically if the server   switches to the other queue while in the middle of a service, then it has to restart that service upon its next visit (see \cite{takagi},\cite{Pinedo}). The jobs are served at both the queues according to   independent and identically distributed  (i.i.d.) service times  and with common rate $\mu$. The customers leave the system only after their service is completed. We define the load factors $\rho_i := \lambda_i / \mu$ $\forall i$ and \textit{assume $\rho:=\rho_1+\rho_2<1$.} Upon visiting a queue, the server attends it until the switching instance  governed by a given two-phase switching  policy described in Section \ref{sec_policy}. 
   \textit{The switching times $\{\Walk_{ik}\}_{k\ge1}  $  from $Q_i$ to $Q_j$   at  various switching instances indexed by $k$  are assumed to be  i.i.d. random variables with mean $\walk_i$; these  are independent of other (previous) events.}

In this work, we aim to identify a sufficiently generic class of stable switching policies for  polling systems. Our further aim is to identify a big class that can be parametrized compactly (using a small set of parameters). For example, a  class general enough to obtain an appropriate and approximate achievable region of the polling system (more details in Section \ref{sec_achievable_region}) or a class general enough to produce required collection of robust policies (see Section \ref{sec_kost_policy}).

 \subsection{A Family of Two-phase Switching Policies, ${\cal T}$}\label{sec_policy}

We consider stationary two-phase switching policies that are employed at each (customer) arrival and departure epoch. \textit{Any such policy depends upon the system state,   the number of waiting customers (say $N_1(t)$, $N_2(t)$)  in the two queues,  at the decision epochs}. While visiting (say) $Q_i$,  the policies are defined using two affine linear functions $\Pi_{iB}(\cdot)$ and $\Pi_{iC}(\cdot)$,

\vspace*{-0.3cm}
\begin{small}
      \begin{eqnarray}
    \Pi_{iB}(N_{i}(t),N_j(t)):=\alpha_i^B N_i(t)+N_j(t)+ \mu \beta_i^B  \label{eqn_fun_pi1}\\
    \Pi_{iC}(N_{i}(t),N_j(t)):= N_i(t)+ \alpha_i^C N_j(t)+\mu \beta_i^C,\label{eqn_fun_pi2}
\end{eqnarray}
\end{small}%
 which respectively determine the beginning and the concluding phase of the server-visit time.  
\textit{Here $\{\alpha_i^m, \beta_i^m\}_{m \in \{B,C\}}$ for both $ i$    are non-positive real constants and represent\footnote{The scaling of $\beta_i^C$, $\beta_i^B$ factors with $\mu$  helps in comparative  analysis across different $\mu$ systems,  however this is not  meant for asymptotic analysis. } the parameters of the policy.} 
Basically the parameters $\{\alpha_i^m, \beta_i^m\}_{m \in \{B,C\}}$ for any $i$ represent the priorities given to the visiting queue ($Q_i$)---the higher values of $B$ parameters  and/ or lower values of $C$ parameters for $Q_i$ imply longer visit times and hence a higher priority to  $Q_i$. 

Now, say the server arrives to  $Q_i$  at time epoch $\Phi$.  During  the subsequent $Q_i$-visit time,  the server only attends to $Q_i$ and  hence     $N_j(t)$  keeps increasing with each arrival to $Q_j$.  On the other hand, the number $N_i(t)$ can increase/decrease with time, but will decrease eventually (as $\rho_i < 1$). These fluctuations eventually determine the end of server $Q_i$-visit phase;   the precise details of the two-phase switching policy are in the following:


 \textbf{Beginning Phase:} The server continues serving the  $Q_i$ as long as the function $\Pi_{iB}(N_{i}(t),N_{j}(t)) < 0$. Once $\Pi_{iB}(N_{i}(t),N_{j}(t))$ becomes zero or positive, the beginning phase, characterized by $\Pi_{iB}(\cdot)$ ends. Hence, the length of the beginning phase at $Q_i$, represented by $B_i$, equals 
 
 \vspace*{-0.3cm}
 \begin{small}
        \begin{eqnarray}
     B_i := \inf \{t \ge \Phi: \Pi_{iB}(N_{i}(t),N_{j}(t)) \ge 0  \}-\Phi.\label{beginning}
    \end{eqnarray}
 \end{small}
 
 \textbf{Concluding Phase:} The concluding phase starts immediately after $B_i$, the server  continues serving $Q_i$ as long as $\Pi_{iC}(N_{i}(t),N_{j}(t))>0$. Once $\Pi_{iC}(N_{i}(t),N_{j}(t))$ becomes zero or negative, the server switches to (immediately starts moving towards)   $Q_{j}$ and the concluding phase at $Q_i$ ends. The length of the concluding phase at $Q_i$ is represented by~$C_i$,

\vspace*{-0.3cm}
 \begin{small}\begin{eqnarray}
        C_i :=  \inf \{t \ge \Phi+B_i : \Pi_{iC}(N_{i}(t),N_{j}(t)) \le  0 \} - (\Phi+B_i).\label{concluding}
    \end{eqnarray}
 \end{small}

In any cycle (say $k$-th cycle), 
the total  visit time of the server at $Q_1$ equals  $M_{1k} := B_{1k} + C_{1k}$. The server leaves $Q_1$ after this visit time $M_{1k}$,  arrives to $Q_2$ further after the random switching time $S_{1k}$ and then  spends  time $M_{2k}$ at $Q_2$   (again  employing a similar switching policy) now at $Q_2$. And this continues.

Note that the functions $\Pi_{iB}(\cdot)$ and $\Pi_{iC}(\cdot)$ are piece-wise constant and their value changes only     at an arrival or a departure epoch. \textit{Thus the arrival and departure epochs form the potential switching epochs.} 
Let ${\cal T}$ represent the collection of two-phase  policies, each determined by   eight parameters:

\vspace{-3mm}
{\small$$
{\cal T} := \left\{ \policy = (\alpha_i^m,\beta_i^m): \alpha_i^m , \beta_i^m \le 0, \  i\in \{1,2\}, \ m \in \{B,C\}     \right\}.
$$}
We will observe that the family ${\cal T}$ is fairly general and includes various policies considered in the literature; some interesting examples are discussed in section \ref{Interesting SP}; for now, observe   that the\textit{ proposed class  can include one phase policies  with only $C$ phase being active, by setting   $\alpha_i^B=0 = \beta_i^B$ for both $i$.}   We next describe the arrival and service process for the polling system and the required assumptions.

\subsection{Arrival and Uninterrupted Service Processes} 

The arrival process to any queue (say $Q_i$) is modeled using a Poisson arrival process with rate $\lambda_i$.  
Let $\Lambda_i (t)$ represent   the  number of arrivals to $Q_i$ until time $t$. 

When the server  attends a queue continuously (without interruptions), \textit{the departure process is referred to as  uninterrupted service process,} which is a part of a renewal process. The uninterrupted service process at $Q_i$ during the $k$-th visit of  the server  is modeled as renewal process\footnote{For mathematical convenience the uninterrupted service process is defined using an infinite sequence, but observe only the first $v$ (i.e., $n\le v$) of them are used if $v$ customers are served during the $k$-th visit; such a definition does not alter anything as the job sizes are i.i.d. across cycles also.} defined via the  i.i.d. service-times $\{ J_{ik,n} \}_{n\ge1} $ (observe here $E[J_{ik,n}] = \nicefrac{1}{\mu}$ for each $i,k,n$). 
The counting process representing the number of departures in $Q_i$ during the $k$-visit of  the server   is represented by  $\Gamma_{ik}  (\cdot)$ and is given by

\vspace*{-0.3cm}
\begin{small}
\begin{eqnarray}
    \Gamma_{ik}  (t) = \sup  \left \{ n :   \sum_{v=1}^n J_{ik, v} \le t  \right \} \mbox{ for any } t \ge 0.
\end{eqnarray}
\end{small}
Observe here that, 
in each visit of the   server to any queue, a new uninterrupted service process begins that  lasts till the switching epoch.

\section{Some Interesting Classes of Policies}\label{Interesting SP}
\subsection{Achievable Region of \cite{WiOpt} and Pareto Frontier}\label{sec_achievable_region}

In \cite{WiOpt}, in a quest to identify the achievable region of the polling system,  defined in terms of customer-average waiting times (or equivalently in terms of stationary  expected number of customers in individual queues by Little's law), the authors analyzed a fluid model under a special class of priority schedulers. We will show that the two-phase switching policies ${\cal T}$ include this class. In section \ref{sec_NI}, we discuss the approximate achievable region using this  class and also  illustrate numerically that the corresponding Pareto frontier converges to the one  derived in \cite[Theorem 8]{WiOpt}. 

The achievable region in terms of customer-average waiting times  is defined as   (see \cite{WiOpt}):
$$
\mathcal{A} =\{ (\bar{w}_1 (\policy), \bar{w}_2 (\policy) ): \textit{ $\policy$ is any stationary  policy}\}, $$
where ${\bar w}_i(\policy)$ is the customer-average of the waiting times (if exists) of the customers of   $Q_i$  under switching policy $\policy$. 
A {\it class ${\cal B}$ of  switching policies is said to be complete,}   if every pair $({\bar w}_1, {\bar w}_2) \in {\cal A}$ is achieved by a switching policy $\policy \in {\cal B}$, i.e., if, \vspace{-4mm}
\begin{eqnarray*}
    {\cal A} &=& \{(\bar{w}_1(\policy), \bar{w}_2(\policy)): \policy \in {\cal B} \},
\end{eqnarray*}
%
%
If such an achievable region is  shown to exist and if a complete class ${\cal B}$ is identified, any optimization problem related to the polling systems can be solved using the feasible performance pairs of  $\mathcal{A}$. For example, the authors in \cite{WiOpt} consider a constrained optimization  as below:
$$ \mbox{Optimize}_{(\bar{w}_1, \bar{w}_2) \in {\cal A}} \,\, \mathcal{F}(\bar{w}_1, \bar{w}_2) \mbox{  s.t. } \bar{w}_1 \leq \upsilon_1,
$$
and observe that it is equivalent to  
$$ \mbox{Optimize}_{ \policy \in {\cal B}} \,\, \mathcal{F}(\bar{w}_1 (\policy), \bar{w}_2 (\policy)) \mbox{  s.t. } \bar{w}_1(\policy) \leq \upsilon_1.$$

In \cite{WiOpt} authors identified one such class of complete schedulers at an appropriate fluid regime. These policies are parameterized by constants, say $\alpha_1, \alpha_2$. Here, the server  switches from $Q_i$ to $Q_j$, at the first time  instance $t$ (after its visit to $Q_i$), when $N_j(t)>N_i(t)\alpha_i$. 

Such a switching policy can be described  by our two-phase policies upon considering  
$$
\alpha^C_i=-\frac{1}{\alpha_i} ,  \ \alpha^B_i = \beta^B_i = \beta^C_i = 0 \mbox{ for each } i, 
$$ 
in \eqref{eqn_fun_pi1}  and \eqref{eqn_fun_pi2}. Hence, for these policies, the length of the beginning phase $B_i$ equals zero. Authors in \cite{WiOpt} also consider exhaustive $(ex)$ policies where the server visiting  $Q_i$ switches only when the queue becomes empty and  $Q_j$ is non-empty. This policy is attained by two-phase policies when all the parameters are zero in \eqref{eqn_fun_pi1} and \eqref{eqn_fun_pi2}. 
The beginning phase equals zero even in case of exhaustive policy, see \eqref{beginning}.  
One can also have exhaustive policy at one queue and a policy with non-zero parameters at other queue (such policies are shown to obtain weak Pareto frontier in \cite[Theorem 6]{WiOpt} and these  switching policies are discussed in section \ref{sec_NI});   one can achieve $\policyex=(ex, \alpha)$ policy, for any $\alpha$,  by setting  in \eqref{eqn_fun_pi1}-\eqref{eqn_fun_pi2}:

\vspace{-3mm}
{\small$$
\alpha_{1}^B = \alpha_1^C = \beta_1^B = \beta_1^C = 0 \mbox{ and } \alpha_2^C = -\frac{1}{\alpha}, \alpha_2^B = \beta_2^B = \beta_2^C = 0. 
$$}

\subsection{Robust Policies \cite{konstantin,Perel2018} }\label{sec_kost_policy}

The focus of the previous class of the switching policies is on expected or customer-average performance.  However, many a times, the focus is on robust design.  
For example, once a queue is above a prescribed (danger) level, one may require to switch immediately and serve that queue.


Authors in  \cite{konstantin,Perel2018} consider switching policies that can be useful in this context. 
They consider a threshold policy, where the server in $Q_i$ considers switching to $Q_j$ once the queue length at $Q_j$ exceeds a threshold $\beta_j$, and waits till the queue length in $Q_i$ falls below a threshold $\beta_i$, and then switches. 
Such policies  can be captured by our two-phase policies   when $\alpha_{i}^B=\alpha^C_{i} = 0$, $\beta_{i}^B=-\beta_j^B$ and $\beta_{i}^C=-\beta_i^C$ for all $i$. Basically, the functions defining these two-phase policies are now given by: 
  \begin{eqnarray}
  \label{Eqn_Robust_policy}
    \Pi_{iB}(N_{i}(t),N_j(t))&=&N_j(t)-\mu \beta_j^B,\nonumber\\
    \Pi_{iC}(N_{i}(t),N_j(t))&=& N_i(t)- \mu \beta_i^C.
\end{eqnarray}
Hence while the sever is visiting $Q_i$, the first phase ends when the queue length at the unattended $Q_j$ is above $\mu\beta_j^B$ and the second phase ends (and the server switches) when the queue length at  $Q_i$ is below $\mu\beta_i^C$.   There is a possibility that the phase $C$ ends before phase $B$; however it is guaranteed that the queue length $N_i$ is at maximum $\mu \beta_i^C$ at all the switching epochs from $Q_i$.

Thus this class of policies    allow the server to immediately jump to the other queue, of course as long as the number of customers in the visiting queue is  below a manageable threshold; hence   they  exhibit  robustness. 

Authors in \cite{konstantin} actually consider a sub-class of policies defined above, where $\beta_i^C = \beta_i^B$ for both $i$; hence their policies are parametrized by two parameters while ours are parametrized by four parameters, but offer more generality. 


\section{Random Systems}\label{approach}

We analyze the polling systems for 
any given  arrival  and departure rates $(\lambda_1, \lambda_2, \mu)$. \textit{The particular focus is on deriving the stability conditions and thereby to  identify the sub-class among ${\cal T}$ that renders the polling system stable for   given  set of rates $(\lambda_1, \lambda_2, \mu)$ and the switching times}.

Our discussions will also include comparison of polling systems with different values of $\mu \in (0, \infty)$. Towards this second purpose, we consider  one  system for each  $\mu$,  while  the load factors  $\rho_i,~i=1,2$ are kept fixed.  Henceforth,  we represent all the quantities with superscript $\mu$ to explicitly show this dependency and call it $\mu$-system. 
We represent the arrival process $\Lambda_i^{\lambda_i} (\cdot) = \Lambda_i^{\lambda_i(\mu)} (\cdot)  = \Lambda_i^{\rho_i \mu} (\cdot)$ at $Q_i$ by $\Lambda_i^\mu (\cdot)$ and the uninterrupted service process for the $k$-th visit of the server\footnote{Note that due to notational complexity, in the remaining part of this section, we use the notation $\Gamma_i^\mu (.)$ instead of $\Gamma_{ik}^\mu (t)$ and correspondingly simplify the notation for service-times.}  by $\Gamma_{ik}^\mu (\cdot)$. 
Thus, for any~$\mu$, we have the following almost surely,

\vspace{-3mm}
{\small\begin{eqnarray}
\label{Eqn_Lam_Gam}
    E\left [ \left . \frac{\Lambda_i^\mu (T)}{\mu} \right | T  \right ] = \rho_i T,  \  \  E \left [ \left . \frac {\Gamma_i^\mu (T) }{\mu} \right | T  \right ] = T  \mbox{  for all } i,   
\end{eqnarray}}%
for any random time $T$ independent of the stochastic processes $\Lambda^\mu(\cdot)$ and $\Gamma^\mu(\cdot)$.  

\hide{
  Let $\{\IA^\mu_{i,n}\}_n$ and $\{J^\mu_{ik,n}\}_n$  represent the exponential inter-arrival times corresponding to the Poisson arrival process $\Lambda^\mu_i(\cdot)$ and inter-departure times (or the job service times) corresponding to the renewal (uninterrupted) service process $\Gamma^\mu_{ik}(\cdot)$, respectively. By definition of exponential random variables we have  $
P(\IA^\mu_{i,n} > 0) = 1,
\mbox{ for all } \mu,i, n
$. For service times,
we assume $P(J^\mu_{ik,n} > 0) = 1
\mbox{ for all } \mu,i, n, k
$. Hence, at any time instance, there can be only one arrival or departure. That is, there are no simultaneous arrivals and departures. Also, we assume that none of the arrival and departure instances coincide.}

\subsection{Evolution of the Polling system} We first describe the evolution of the polling system with service rate $\mu$ and introduce some of the notations to present the main results of this work.  
Without loss of generality, we define \textit{a cycle to be the time period between subsequent server-visit epochs of $Q_1$}. We represent the duration of $k$-th cycle by $\Psi^\mu_k$. Let $\Phi^\mu_k := \sum_{n < k} \Psi^\mu_n $ represent the starting epoch of $k$-th cycle. Recall  $N^\mu_i(t)$ represents the number of customers in  $Q_i$ at time $t$ and let $N_{ik}^\mu := N^\mu_i(\Phi_k)$ represent this  number   at $\Phi^\mu_k$.  \textit{Let  $S_{ik}^\mu$ be the random switching time from $Q_i$ to $Q_{j}$ in the $k$-th cycle with mean $s_{i}^\mu$}. Similarly, let  $M_{ik}^\mu, B_{ik}^\mu$, and $C_{ik}^\mu$ be the total visit time, and the durations of the beginning and the concluding phases of the visit time, respectively. Hence, $M^\mu_{ik}=B^\mu_{ik}+C^\mu_{ik}$. \hide{Recall,   affine functions 
\begin{eqnarray}
    \Pi_{iB}(t) &:=&\alpha_i^B N_i^\mu(t)+N_j^\mu(t)+\mu \beta_i^B, \label{beginning1}\qquad\mbox{ and} \\ \Pi_{iC}(t) &:=& N_i^\mu(t)+\alpha_i^C N_j^\mu(t)+\mu \beta_i^C.\label{concluding1}
\end{eqnarray} 
define the switching instances.} 

\begin{figure}[htb!]\centering
\includegraphics[scale=0.37]{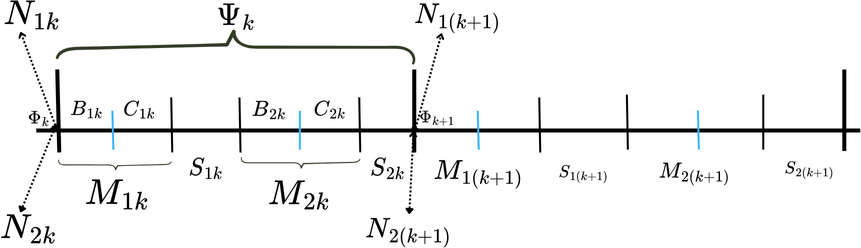}
\caption{Polling system with two-phase switching policy}\label{Two queue polling system}
\end{figure}

With\textit{ general i.i.d. service and switching times, it is clear that our process is not Markovian. However  our switching policies are Markovian: basically the decisions depend upon the   `state' at the decision epoch (say~$t$)},  $ (N_1^\mu (t), N_2^\mu (t), P^\mu(t) ) $ with $P^\mu(t)$ taking values in the set $\{(1, B_1), (1, C_1), S_1, (2, B_2), (2, C_2), S_2 \}$ that represents various  phases of the system  (see Fig. \ref{Two queue polling system}); here $(i,B_i)$  and  $(i,C_i)$ respectively represent the beginning and concluding phases in $Q_i$ while $S_i$ represents the switching phase from~$Q_i$. 
\hide{
, we define six phases of our polling system during each cycle.  depicts these six phases: (i) beginning phase in $Q_1$ described by $(1, B_1)$; (ii) concluding phase in $Q_1$ described by $(1, C_1)$; (iii) switching phase from $Q_1$ to $Q_2$ described by $S_1$; (iv) beginning phase in $Q_2$ described by $(2, B_2)$;(v) concluding phase in $Q_2$ described by $(2, C_2)$; (vi) switching phase from $Q_2$ to $Q_1$ described by $S_2$. 
}
\vspace{-0.2cm}

 We now present the evolution of the number of customers in any queue. 
For ease of notation, define 
\begin{equation}
\small D^\mu_{1k} :=  \Phi^\mu_k + M^\mu_{1k}, \mbox{ and, }  
D^\mu_{2k}:= \Phi^\mu_k+ M^\mu_{1k} +S^\mu_{1k}+M^\mu_{2k} ,  \hspace{3mm}\label{Eqn_Dik}
\end{equation} 
 respectively as the departure epochs from $Q_1$ and $Q_2$ in $k$-th cycle. Then the number  of customers in $Q_1$ evolve as:

\vspace{-2mm}
{\small\begin{eqnarray}
 N^\mu_{1}(t) = \left\{ \begin{array}{ll}
      {N^\mu_{1k} + \Lambda_1^{\mu}(t) - \Lambda_1^{\mu}(\Phi^\mu_k)  - \Gamma^\mu_{1k} (t-\Phi_k^\mu)} \vspace{0.1cm}\\ &\hspace{-13mm}
    \mbox{ if }   \Phi^\mu_k \le t \le D^\mu_{1k}, \vspace{0.21cm}\\
        {N_1^\mu(D^\mu_{1k}) + \Lambda_1^{\mu}(t) - \Lambda_1^{\mu}(D^\mu_{1k})}\\  & \hspace{-14mm}
    \mbox{ if }  D^\mu_{1k} <t< \Phi^\mu_{k+1},   
 \end{array}\right. \hspace{-10mm}\nonumber \\
 \mbox{for some }  k.
 \label{Queue1} 
\end{eqnarray}}
In the above, $\Lambda_1^\mu(t) - \Lambda_1^\mu(s)$ represents the number of arrivals in $Q_1$ in the time duration $(t-s)$ for $s \le t$.  Since the uninterrupted service process is assumed to be non-resume preemptive, $\Gamma^\mu_{1k}(t-s)$ equals the number of departures from $Q_1$ in time duration $(t-s)$; recall we have one such process for each $k$ which ends at $D_{1k}^\mu$. Once the server departs from $Q_1$ at epoch $D_{1k}^\mu$, there are only arrivals at $Q_1$ (and no departures). Similarly, the number of customers in $Q_2$ evolve as below:

\vspace{-3mm}
\begin{small}
    \begin{eqnarray}
    N^\mu_{2}(t) = \left\{ \begin{array}{ll}
      {N^\mu_{2k}+ \Lambda_2^{\mu}(t)  - \Lambda_2^{\mu}(\Phi_k^\mu)}    
     & \hspace{-11mm}
    \mbox{if  }   \Phi^\mu_k \le t \le D^{\mu}_{1k} +S^{\mu}_{1k},  \vspace{0.21cm}\\
\substack{N_{2}^\mu(D^{\mu}_{1k} +S^{\mu}_{1k})+ \Lambda_2^{\mu}(t)-\Lambda_2^{\mu}(D^{\mu}_{1k} +S^{\mu}_{1k})\\ - \Gamma^\mu_{2k} (t-D^{\mu}_{1k} -S^{\mu}_{1k})} \\ &\hspace{-11mm} 
    \mbox{if }  D^{\mu}_{1k} +S^{\mu}_{1k} < t \le D^{\mu}_{2k},  \vspace{0.21cm}\\
   \hspace{-2mm}    {N_2(D^{\mu}_{2k}) + \Lambda_2^{\mu}(t) - \Lambda_2^{\mu}(D^{\mu}_{2k})} \\& \hspace{-11mm} 
    \mbox{if  }D^{\mu}_{2k}<t< \Phi^\mu_{k+1},  
 \end{array}\right. \hspace{-10mm} \nonumber \\\label{Queue2}
 \mbox{ for some } k. 
\end{eqnarray}
\end{small}%
In the above, $\Lambda_2^\mu(t) - \Lambda_2^\mu(s)$ and  $\Gamma^\mu_{2k}(t-s)$  represent the number of arrivals and departures in $Q_2$  in duration of length  $t-s$, as before.  Note that we have departures in $Q_2$  only during   $ D^{\mu}_{1k} +S^{\mu}_{1k} < t \le D^{\mu}_{2k}$, i.e., during its server-visit phase.

\subsection{Markov Chain at Server Visit Epochs} 
With  Poisson arrival process, $\{(N^\mu_{1k},N^\mu_{2k})\}_{k\in \mathbb{N}}$ is a discrete-time Markov chain for every $\mu$ and so is the case with the normalized chain, $\{{\bf X}_k^\mu\}_{k\in \mathbb{N}}$, where  

\vspace{-4mm}
\begin{small}
    $$
{\bf X}_k^\mu:=(\Theta^\mu_{1k},\Theta^\mu_{2k}), \mbox{  }\Theta^\mu_{ik}= \Theta_i^\mu(\Phi_k),  \ \Theta_i^\mu (t) := \frac{N^\mu_{i}(t)}{\mu},i=1,2.$$ 
\end{small}%
We refer to $\{{\bf X}_k^\mu\}_{k\in \mathbb{N}}$  as  the \textit{Palm Markov chain}. In this section, we aim to derive the stability of this chain and then eventually that of the joint stochastic process ${\bf X}^\mu  :=  ({\bf X}^\mu (t) )_{t\ge 0}$, defined for all~$t$, where  ${\bf X}^\mu(t) := (\Theta_1^\mu(t), \Theta_2^\mu(t))$. 

\hide{Therefore, we find the condition for which this Markov chain will be stable. Now $\{{\bf X}_k^\mu\}_{k\in \mathbb{N}}$ is Palm Markov chain of the original stochastic process $(N_1^\mu(t),N_2^\mu(t))$. Thus if we get the stability of the Palm Markov chain we can also get the stability of original stochastic process.

Before going to stability, let us discuss some important results for our $\cal T$ switching policy, which will help to find the stability conditions precisely. For the stable polling system, as the visiting phase $M_{ik}^\mu$ is always positive for both queue, either of the beginning phase $B_{ik}^\mu$ and the concluding phase $C_{ik}^\mu$ is positive or both of them are positive. Hence, we have following lemma.
\begin{lem}
     If both the beginning and concluding period are positive for our $\cal T$ switching policy, i.e.,  $M^\mu_{ik}> B^\mu_{ik} >0$, for all $k=1,2,\ldots$ and $i=1,2$, then $\alpha_{i}^B\alpha_{i}^C > 1, i=1,2$.\\
 \end{lem}
 }

 From \eqref{concluding}, 
the affine functions $\Pi_{iC}$ with $i=1,2$ determine the concluding phases  $C_{ik}^\mu$ of the $k$-th visit time. Let random variables,  $\mathcal{H}_{1k}^\mu$ and $\mathcal{H}_{2k}^\mu$ respectively represent the  values of this function at the switching instances in the two queues. Using the departure times of \eqref{Eqn_Dik}
 we have
 
 \vspace*{-0.4cm}
 \begin{small}
         \begin{eqnarray}
       \Pi_{1C}(N^\mu_{1}(D^\mu_{1k}),N^\mu_{2}(D^\mu_{1k}))&&\nonumber\\
       &&\hspace*{-2.50cm}= N^\mu_{1}(D^\mu_{1k})+\alpha_{1}^C N^\mu_{2}(D^\mu_{1k}) +\mu \beta_{1}^C=\mathcal{H}_{1k}^\mu,\label{visiting1}
       \\ 
     \Pi_{2C}(N^\mu_{2}(D^\mu_{2k}),N^\mu_{1}(D^\mu_{2k}))&&\nonumber\\&&\hspace*{-2.5cm}= N^\mu_{2}(D^\mu_{2k})+ \alpha_{2}^C N^\mu_{1}(D^\mu_{2k})+\mu \beta_{2}^C=\mathcal{H}_{2k}^\mu. \hspace{4mm}\label{visiting2}
\end{eqnarray}
 \end{small}
We now assume the following about these residual values:
\begin{itemize}
      \item [{\bf R.1:}] Assume 
    $E[\mathcal{H}_{ik}^\mu],$
      is bounded uniformly in $k$, for each $i.$ 
 \end{itemize}
 
If the concluding $C$ phase is strictly positive for both the queues at all the departure epochs, \textit{then the above is trivially true}, as explained in the immediate next. Let $\Pi_{iC} (t^+)$ and $\Pi_{iC} (t^-)$ respectively represent the value of $\Pi_{iC}$ function in \eqref{eqn_fun_pi2}  immediately after and  before the departure epoch $D^\mu_{ik}$. Observe that for strictly positive concluding phases,  by definition, $\Pi_{iC} (t^-) < 0$ and $\Pi_{iC} (t^+) \ge 0$  and such a transition is triggered either by an arrival or by a departure. Thus in fact $\mathcal{H}_{ik}^\mu $ are bounded almost surely by $|\alpha_i^C|+1$, and hence, {\bf R.1} is trivially true.  

Therefore {\bf R.1} is true for  the class of policies considered in Section \ref{sec_achievable_region}, the ones that include the Pareto complete class at fluid limit. 

From \eqref{Eqn_Robust_policy}, it is clear that  $E[\mathcal{H}_{ik}^\mu] \le \mu|\beta_i^C |$ for all $k$ and thus {\bf R.1} is also satisfied for the robust class described in Section \ref{sec_kost_policy}.

The assumption {\bf R.1} can be a restriction only for the subclass of policies for which the concluding phases in some cycles are zero (possible if $C$ phase ends before $B$ phase sometimes); Section \ref{sec_kost_policy} describes one such class, but the assumption {\bf R.1} is satisfied for this class  directly using the uniform boundedness of the visiting queue at the server-departure epochs. It is also intuitive that one cannot anticipate stability without   assumption {\bf R.1}.  

Now to determine the stability of the Palm Markov chain $\{{\bf X}_k^\mu\}_{k\in \mathbb{N}}$,
we use the Foster stability criteria by constructing an appropriate Lyapunov function (see e.g., \cite{Lyapunov}). In particular, we consider the following as Lyapunov function
\begin{eqnarray}
  V({\bf X}_k^\mu)=\rho_2\Theta^\mu_{1k}+(1-\rho_1)\Theta^\mu_{2k}. \label{Eqn_Lyapunov}  
\end{eqnarray}
Then, the mean drift is given by 

\vspace*{-0.3cm}
\begin{small}
    \begin{eqnarray}
    \Delta V({\bf x}^\mu)&=& E[V({\bf X}^\mu_{k+1})-V({\bf X}^\mu_{k})| {\bf X}^\mu_{k}={\bf x}^\mu]\nonumber\\
    &=& E[\rho_2(\Theta^\mu_{1(k+1)}-\Theta^\mu_{1k})\nonumber\\&&\hspace{-2mm}+(1-\rho_1)(\Theta^\mu_{2(k+1)}-\Theta^\mu_{2k})| {\bf X}^\mu_{k}={\bf x}^\mu], \hspace{5mm}\label{Froster}
\end{eqnarray} 
\end{small}
where ${\bf x}^\mu:=(\theta_1^\mu,\theta_2^\mu)$. 

Now at any $k$-th cycle, we have the following relation for the number of customers at two consecutive cycles in $Q_i$ of the random system (see \eqref{Queue1} and \eqref{Queue2}):

\vspace*{-0.5cm}
\begin{small}
\begin{eqnarray}
    N^{\mu}_{i(k+1)}   &=&  N_{ik}^{\mu}    + \Lambda_{i}^{\mu} \left(M^\mu_{1k}+S^\mu_{1k}+M^\mu_{2k}+S^\mu_{2k}\right)\nonumber\\&&  - \Gamma_{i}^\mu \left(M^\mu_{ik}\right),~i=1,2. \label{THETAi}
    \end{eqnarray}
\end{small}

   Hence, by normalizing the above equations for $i=1,2$ with $\mu$
   and  taking the expectations $E[.|{\bf X}^\mu_{k}={\bf x}^\mu]$, we obtain the following, using \eqref{Eqn_Lam_Gam}, 

   \vspace*{-0.5cm}
\begin{small}
     \begin{eqnarray*}
   E[\Theta_{i(k+1)}^\mu|{\bf X}^\mu_{k}={\bf x}^\mu] \nonumber\\&& \hspace{-25mm}=\theta_{i}^\mu+\rho_i\left(E[M^\mu_{1k}+M^\mu_{2k}|{\bf X}^\mu_{k}={\bf x}^\mu]+s_1^\mu+ s_2^\mu\right)\\&&-E[M^\mu_{ik}|{\bf X}^\mu_{k}={\bf x}^\mu].\label{Implicit1}
 \end{eqnarray*}
    \end{small}%
Consequently, from \eqref{Froster}, we have (recall $\rho =\rho_1+\rho_2$)

\vspace*{-0.3cm}
\begin{small}
  \begin{equation}
       \Delta V({\bf x}^\mu)=(\rho-1)E[M_{2k}^\mu|{\bf X}^\mu_{k}={\bf x}^\mu]+\rho_2(s_1^\mu+ s_2^\mu). \hspace{3mm}\label{Froster1}
  \end{equation}
  \end{small}
  We now provide the conditions under which  the  two-phase switching policy $\policy \in {\cal T}$ renders the Palm chain stable.
\begin{thm}\label{thm_palmstability}
{\bf [Palm stability]}
Assume {\bf R.1}. Then the sufficient condition for stability of the Markov chain   $\{{\bf X}_k^\mu\}_{k\in \mathbb{N}}$ with switching policy $\policy \in {\cal T}$  is that \vspace*{-0.1cm}
$$
\rho_1+\rho_2 < 1 \mbox{ \ and \ } \alpha_1^C \alpha_2^C < 1.
$$
\vspace*{-0.3cm}\end{thm}
\begin{proof}
    The proof is in Appendix.
\end{proof}
The intuition behind taking such a Lyapunov function is the following. Any $k$-th cycle, by our definition,  starts with server-visit of  $Q_1$ and length $\Psi^\mu_k$ of the cycle is the time lapsed  between two consecutive visits to $Q_1$.  The total visit time at $Q_2$ within this $\Psi^\mu_k$,  is given by $M_{2k}^\mu$. This depends upon the number of customers in both the queues  at the start of $k$-th cycle,  i.e., on $\Theta_{1k}^\mu$  and $\Theta^\mu_{2k}$ and also upon  the visit time to  $Q_1$, given by $M_{1k}^\mu$. Hence, it was sufficient to construct  Lyapunov function  \eqref{Eqn_Lyapunov} in such a way that the drift of the function would be influenced only by the change in the second visit time of the cycle (i.e., visit time of  $Q_2$), see~\eqref{Froster1}.

Let ${\cal T}_s$ represent the  subclass of stable switching policies: \vspace*{-0.2cm}
\begin{equation}
  {\cal T}_s  := \left \{ \policy \in {\cal T} :  \alpha_1^C\alpha_2^C<1 \right \}.
\end{equation}
\textit{Interestingly, once $\rho < 1$, the stability depends only on two parameters $(\alpha_1^C, \alpha_2^C)$   of the policies.}

We  now prove that the first moments related to Palm chain converge as  the cycle number increases. 
\begin{thm}\label{thm_palmmoments}
{\bf [Palm First moments]}
Assume {\bf R.1}, $\rho < 1$. For any $\policy \in {\cal T}_s$, as the cycle number $k \to \infty$, we have:
\vspace*{-0.1cm}
\begin{eqnarray*}
  E[\Theta_{ik}^\mu] &\to& \theta_i^* + o(\nicefrac{1}{\mu}) \mbox{ for each }  i \in {1, 2}, \mbox{ and }  \\
  E[\Psi_k^\mu] &\to& \psi^* < \infty \mbox{ and hence, } \sup_k  E[\Psi_k^\mu]  < {\overline \psi} < \infty,
 \end{eqnarray*}
 \vspace*{-0.3cm}\\
where $\theta_1^*, \theta_2^*$ are defined in \eqref{stabletheta1}, \eqref{stabletheta2}, and $\psi^*$ and ${\overline \psi}$ are some finite  ($\mu$ dependent) constants.  
\end{thm}
\begin{proof}
    The proof is in Appendix.
\end{proof}

\subsection{Time Stationarity}  
The time epochs $R_m^\mu$  are called regeneration epoch  for a stochastic process ${\bf X}^\mu(t)$, if   the following two stochastic processes are stochastically  identical (e.g., \cite{konstantin2023,Anurag}):
$$
{\bf X}^\mu(R_1^\mu+ \cdot) \stackrel{d}
{=} {\bf X}^\mu(R_m^\mu+\cdot) \mbox{ for each } m \ge 1. 
$$ 
Observe that the regeneration times $T_m^\mu := R_{m}^\mu- R_{m-1}^\mu$, for $m >1$,  are i.i.d. by definition. 
The process ${\bf X}^\mu(\cdot)$ is called positive recurrent if $E[T_m^\mu ] < \infty $ and in such  cases (e.g., \cite{Anurag}) the time stationary limit distribution $(\nu_{{\bf x}})$ exists---for each ${\bf x}$ the following is true almost surely:

\vspace*{-0.3cm}
\begin{small}
\begin{eqnarray}
\label{Eqn_stat_dist}
    \lim_{t \to \infty} \frac{\bigints_o^t 1_{\{ {\bf X}^\mu(t) = {\bf x} \}} dt}{t}  = \frac{  E\left [ \bigints_{R_{m-1}^\mu}^{R_m^\mu} 1_{\{ {\bf X}^\mu(t) = {\bf x} \}} dt \right ] }{E[T_m^\mu]} =: \nu_{{\bf x}}.
\end{eqnarray}
\end{small}

Towards proving time stationarity of our polling system ${\bf X}^\mu (t)$, representing  the normalized versions of \eqref{Queue1}-\eqref{Queue2}, we begin with identifying a regeneration epoch  using Theorem~\ref{thm_palmstability}.  

Define the following server visit epochs (to $Q_1$), recursively (assuming that the server starts at $Q_1$, i.e., $\Phi_0^\mu=0$)

\vspace*{-0.2cm}
\begin{small}
    \begin{eqnarray*}
  R_{1}^\mu  &:=& \inf \left \{ \Phi_k^\mu : k \ge 1, \Theta_{ik}^\mu  = \theta^*_i \mbox{ for both } i  \right \}  \mbox{, and, } \\
  R_m^\mu &:=&   \inf \left \{ \Phi_k^\mu > R_{m-1}^\mu :   \Theta_{ik}^\mu  = \theta^*_i \mbox{ for both } i  \right \}, \forall \  m > 1.
\end{eqnarray*}
\end{small}

By Theorem \ref{thm_palmstability}, $(\theta_1^*, \theta_2^*)$ are recurrent states of the stationary Palm chain, and hence, clearly $R_m^\mu$ with $m > 1$,  defines the regeneration epochs for the system, due to  Poisson arrivals (see e.g., \cite{Anurag}).
Define the  $m$-th regeneration cycle $T_m^\mu := R_{m}^\mu - R_{m-1}^\mu$, and  
let $\xi_m^\mu$ represent the `Palm cycle' number accounting to $m$-th regeneration epoch, $R_m^\mu$. Then the  regeneration epoch $R_m^\mu$ and  regeneration cycle $T_m^\mu$  can be rewritten as:
\begin{eqnarray*}
    R_m^\mu = \Phi_{\xi_m^\mu}^\mu \mbox{ and   }  T_m^\mu   = \sum_{k=\xi_{m-1}^\mu}^{\xi_m^\mu} \Psi_k^\mu.
\end{eqnarray*} 

Thus, by Theorems \ref{thm_palmstability}-\ref{thm_palmmoments}, with $\nu^P$ representing the stationary distribution of the Palm chain,  for any $m  > 1$ we have:
 \begin{eqnarray}
    E[ T_m^\mu ] \le \sup_k E[\Psi_k^\mu] E[\xi_{m}^\mu-\xi_{m-1}^\mu] < \frac{{\overline \psi} }{ \nu^P_{\theta_1^*, \theta_2^*}}  < \infty. 
\end{eqnarray} 

We thus have proved the following   using the results of~\cite[Section 3.5.1]{Anurag}:
\begin{thm}\label{timestability} {\bf [Time stationarity]} \textit{Assume {\bf R.1}, $\rho < 1$ and consider any $\policy \in {\cal T}_s$.  Then the  stochastic process ${\bf X}^\mu(t)$ representing the two-queue polling system is positive recurrent  and satisfies 
\eqref{Eqn_stat_dist}.}
\end{thm}

\section{Numerical Illustration}
\label{sec_NI}

In this section we would like to obtain individual-queue-centric performance measures like the stationary expected number of
waiting customers in each (say, $i$-th) queue for any $\policy \in {\cal T}_s$:
 \begin{eqnarray}
  \lim_{K\to \infty }  \frac{1}{K} \int_0^K N_i^\mu (t) dt = E_{\nu_\policy} [N_{i*}^\mu]   \mbox{ (when it exists) }. \label{Eqn_Stat_moments}
\end{eqnarray}
In the above, $N_{i*}^\mu$ denotes the number of customers in $Q_i$ in steady-state and $\nu_\policy$ denotes the stationary distribution under the given stable switching policy $\policy \in {\cal T}_s$.
In the strict sense, one first needs to prove the existence of   stationary moments in \eqref{Eqn_Stat_moments}, before discussing them. 
The precise proof of this result is postponed to journal version, however we will now provide partial theoretic and partial heuristic arguments  towards the conjecture that such moments exists for $\policy \in {\cal T}_s.$ It is not difficult to observe the  following (recall $\Phi^\mu_k$ is the server-visit epoch to~$Q_1$) for any $k$:
    \begin{eqnarray}
    \sup_{ t \in [\Phi_k^\mu, \Phi^\mu_{k+1}] } \Theta_1^\mu (t)   \le   \Theta_1^\mu(\Phi_k^\mu) + \Lambda_1^\mu (\Psi_k^\mu)  \mbox{ almost surely, }
    \label{Eqn_Upperbound}
    \end{eqnarray}
   
     and thus  
    $
    E[\Theta_1^\mu(t)] \le  \sup_{k} \left ( E[\Theta_{ik}^\mu]  + 
  E[\Psi_k^\mu]   \right ) \mbox { for any }  t.
$
By Theorem \ref{thm_palmmoments}, these   moments   are bounded  and so we expect  to prove (recall $T$ is a typical regeneration epoch):

{ $$
E\left [ \int_0^T \Theta_1^\mu (t) dt \right  ] < \infty;
$$}%
if the above is proved, then  the results of regenerative processes confirm the existence of the stationary moments in \eqref{Eqn_Stat_moments} (similarly for $Q_2$ consider   a Palm chain at server-visit epochs to $Q_2$).  For now we assume the conjecture and proceed towards numerical observations.


Now, among the  stable class ${\cal T}_s$, we are interested in finding a \textit{suitable subclass of policies that generates the Pareto frontier for the achievable region (the policies under which the performance vectors are strongly or weakly efficient, in that one can't improve one of them without hampering the other or one cannot strictly improve both of them respectively).} 
By  \cite[Theorem 6]{WiOpt}, the
  Pareto-complete sub-class at fluid limit is given by:
$$
\hspace{-2mm}
{\cal P}_s = \{ \policyex:\policyex=
(ex, \alpha^C) \mbox{ or  }(\alpha^C, ex), \ \  \alpha^C\le 0   \}.
$$
We now examine numerically if   ${\cal P}_s$ sub-class generates Pareto frontier (among ${\cal T}_s$) specially for sufficiently large $\mu$.


   

\hide{
For any $\beta \in {\cal T}_s$ by virtue of Theorem \ref{timestability}, we have
 \begin{eqnarray*}
 E_{\pi_\beta} [N_{i*}^\mu]= \lim_{K\to \infty }  \frac{1}{K} \int_0^K N_i^\mu (t) dt. 
\end{eqnarray*}}
We use \eqref{Eqn_Stat_moments}  to numerically obtain the stationary performance measures using Monte Carlo simulations.
The simulated events are the arrivals and the departures of the customers and the switching of the server from the visiting queue to the other.

\begin{figure}[htb!]\centering
\vspace{-1mm}
\begin{minipage}{0.45\textwidth}
    \includegraphics[height=4cm, width=4.2cm]{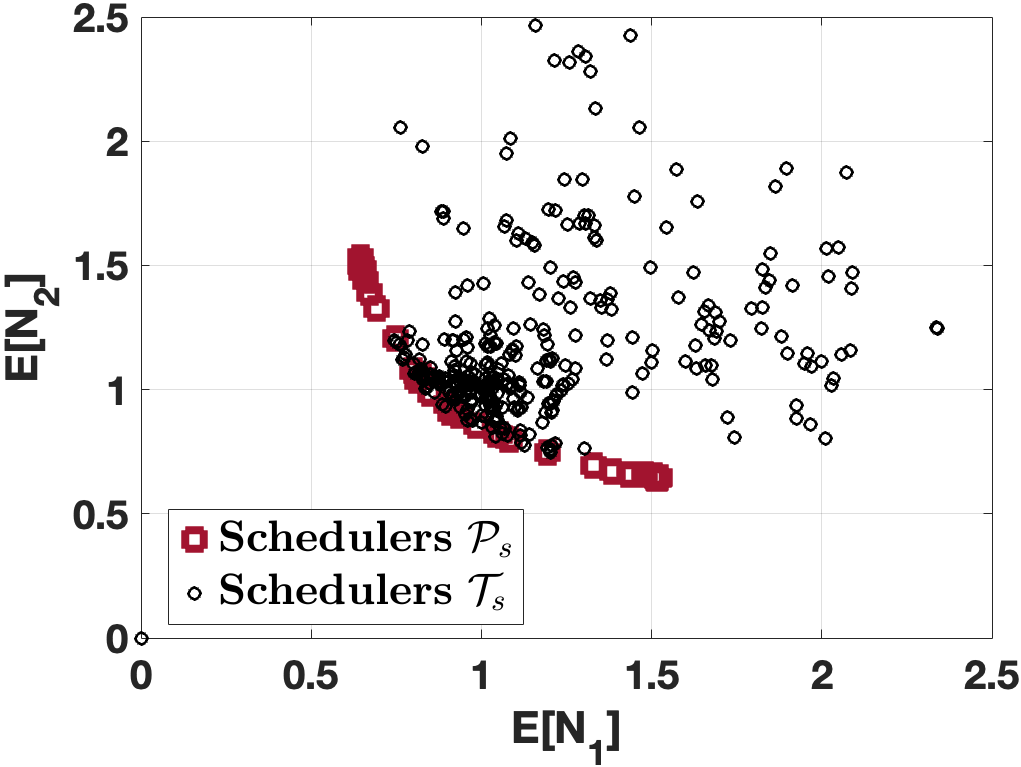}
\end{minipage}%
\begin{minipage}{0.45\textwidth}
\hspace{-36mm}
\includegraphics[scale=0.12]{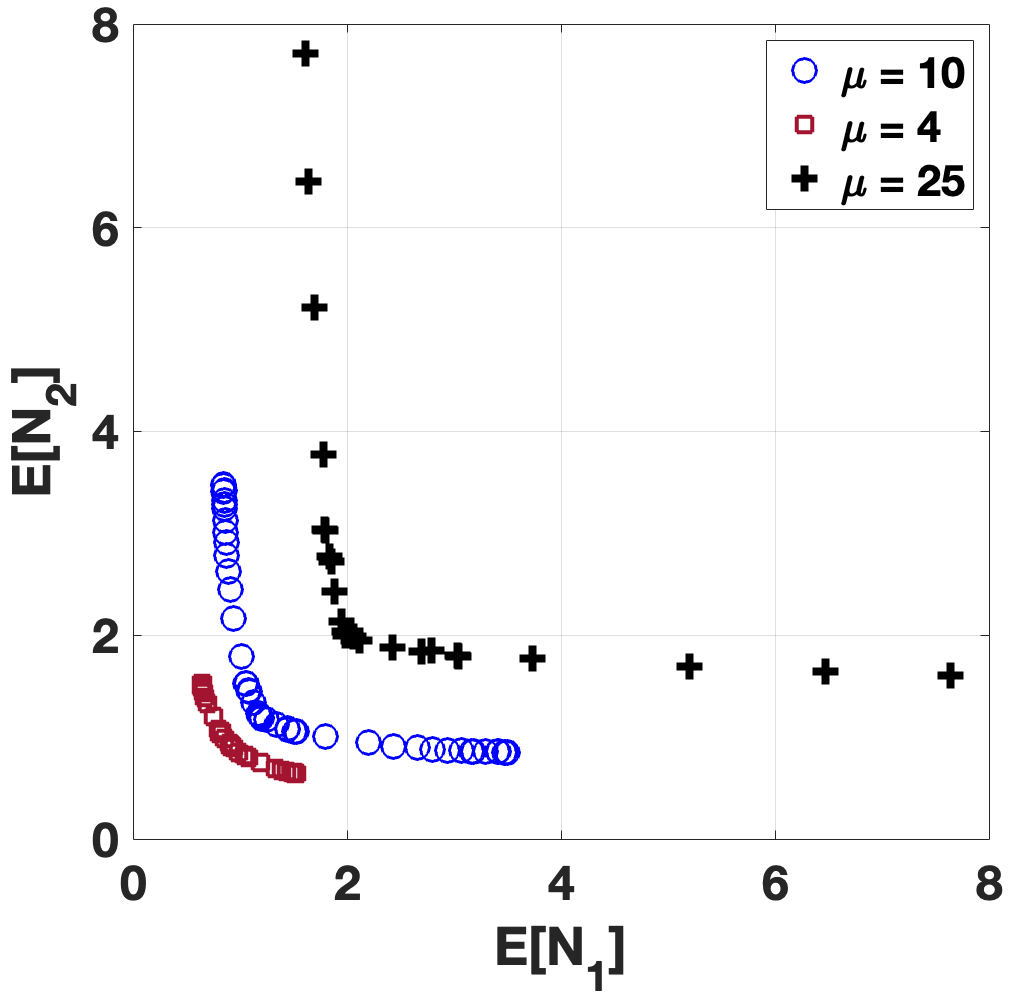}    
\end{minipage}%
\vspace{-1mm}
\caption{{\bf Left figure:}  Performance comparison of Pareto class  ${\cal P}_s$ with   ${\cal T}_s$. \  
{\bf Right figure:} Pareto frontier for different $\mu$.}
\label{fig_Pareto frontier}
\end{figure}
In the left sub-figure of Fig. \ref{fig_Pareto frontier},
we plot the estimates of the stationary expected number of customers in each queue for  some switching policies drawn randomly from class ${\cal T}_s$ (see black small circles). In the same figure and for the same set of parameters, we also plot the performance measures for some policy $\policyex  $  from supposed Pareto class ${\cal P}_s$ (see red squares). As seen from the figure the performance of the sub-class ${\cal P}_s$ is indeed on the Pareto frontier. 

In the right sub-figure of Fig. \ref{fig_Pareto frontier}, we plot the performance measures for some policy $\policyex \in {\cal P}_s$ for different values of $\mu$. This supposedly represents the Pareto frontier for different $\mu$, and as seen from the figure, the Pareto frontier indeed approaches the weak frontier obtained in \cite[Theorem 8]{WiOpt} for the fluid limit. 

{\section{Future directions} 
 One of the most interesting future direction of the current work is to analyze our polling setup for large $\mu$-systems. By Theorem \ref{thm_palmmoments},    
  $
  \lim_{\mu \to \infty} \lim_{k\to\infty }E[\Theta_{ik}^\mu] = \theta_i^* \mbox{ for both } i $. 
  In fact by using the results of \cite{Bertsimas} and examining the small set \eqref{Eqn_small_set} used for proving stability  in Appendix, one can observe that as $\mu \to \infty$  the random sequence $\Theta_{1k}^\mu$ in steady state concentrates near $(\theta_1^*, \theta_2^*)$ for large $\mu$ (recall scaling with $\mu$). These results can probably be used in constructing a  robust achievable region (see \eqref{Eqn_Upperbound}), defined probably in  terms of probabilities of queue-lengths exceeding some given thresholds. Working further on regenerative arguments (second moments), one can directly derive the time-average or stationary performance and can derive their limit in fluid regime, or when $\mu \to \infty.$ We anticipate such limits to have simple closed-form expressions which can provide much clearer insights for systems with large $\mu$.  }

\section{Conclusion}
This work introduces a geneal class of two-phase Markovian switching policies in a two-queue polling system with infinite buffer, general service times and non-zero general switch-over times. We identify the subclass of the policies that are stable; interestingly once the load factor is fixed at a value less than one, the stability condition remains the same irrespective of other  system parameters (like service and arrival rates).

Eight parameters can represent our policies in the most generic form, while the stability of the policy is determined by only two parameters. 
We observe that  the  time limit of some first order moments (related to quantities at server visit epochs) converge  as  the service rate $\mu$ increases to $\infty$.  Using such observations,  we discuss potential connections to different varieties of achievable regions.

\section*{Appendix}
\begin{figure*}
\begin{small}
    \begin{eqnarray}
   a_{11}&=& \frac{(\rho_1+\rho_2-1)-(\rho_1-1)\rho_1\alpha_{2}^C-\rho_1\rho_2\alpha_{1}^C\alpha_{2}^C}{\spio\spit},  \ 
    a_{12} \ = \ \frac{(1-\rho_1)\rho_1-(1-\rho_1)(1-\rho_2)\alpha_{1}^C}{\spio\spit}\nonumber\\
     a_{21}&=&\frac{-\rho_1\rho_2\alpha_{2}^C-(\rho_2-1)\rho_2\alpha_{1}^C\alpha_{2}^C }{\spio\spit},  \  \   
     a_{22} \  = \  \frac{(1-\rho_1)(\rho_2-1)+\rho_2(1-\rho_2)\alpha_{1}^C-(\rho_1+\rho_2-1)\alpha_{1}^C\alpha_{2}^C}{\spio\spit}\nonumber\\
      a_{13}&=&\frac{(1-\rho_1-\rho_2)\beta_{1}^C}{\spio\spit} -\frac{\rho_1\beta_{2}^C+\rho_1s_{1}^\mu}{\spit}-\rho_1s_{2}^\mu,\ 
     \
     a_{23} \ = \ \frac{\alpha_{2}^C(\rho_1+\rho_2-1)\beta_{1}^C}{\spio\spit}  +\frac{(1-\rho_2)\beta_{2}^C+\rho_1s_{1}^\mu\alpha_{2}^C}{\spit}-\rho_2s_{2}^\mu. \label{stabletheta2}
\end{eqnarray}
\end{small}
 \hrule
\end{figure*}

{\bf Proof of Theorem 
\ref{thm_palmstability}:}
\label{Proof_thm_1} 
Let ${\bm \varrho} := (\rho_1,\rho_2).$
We begin with obtaining  the conditional expectation $E[M_{2k}^\mu|{\bf X}^\mu_{k}={\bf x}^\mu]$,  conditioned on   initial state  
 ${\bf x}^\mu = (\theta_1^\mu, \theta_2^\mu)$ 
towards analysing
\eqref{Froster1}.
Taking   conditional expectations $E[.|{\bf X}^\mu_{k}={\bf x}^\mu]$ of 
  equation \eqref{visiting1}    (after normalizing with $\mu$) and rearranging we have (also using    \eqref{Eqn_Lam_Gam}, \eqref{Eqn_Dik}, \eqref{Queue1}, and \eqref{Queue2}):

\vspace{-2mm}
{ \small
    \begin{eqnarray}
E[M^\mu_{1k}|{\bf X}^\mu_{k}={\bf x}^\mu] &=& \frac{\theta_{1}^\mu+\alpha_{1}^C \theta_{2}^\mu+\beta_{1}^C}{\spio}-\frac{E[\mathcal{H}_{1k}^\mu|{\bf X}^\mu_{k}={\bf x}^\mu]}{\mu\spio},\nonumber   \\
\mbox{where, }\spio &:=&  {(1-\rho_1)-\alpha_{1}^C\rho_2},
\label{stopping3} 
 \end{eqnarray}}
 as from \eqref{Queue1}, we (for example) have:

\vspace{-3mm}
{\small$$
E\left [ \left . \frac{N_1(D_{1k}^\mu)}{\mu} \right  | {\bf X}^\mu_{k}={\bf x}^\mu \right ] = \theta_{1}^\mu +  (\rho_1 - 1) E[M^\mu_{1k}|{\bf X}^\mu_{k}={\bf x}^\mu].  
$$}
Observe $\spio > 0$, as $\alpha_1^C \le 0$. 
Now considering  equation \eqref{visiting2}, proceeding as above  and using \eqref{stopping3},  we have

\vspace{-2mm}
\begin{small}
\begin{eqnarray}
 E[M^\mu_{2k}|{\bf X}^\mu_{k}={\bf x}^\mu]&&\nonumber\\&& \hspace*{-2.95cm}= \frac{(1-\alpha_{1}^C\alpha_{2}^C)\left((1-\rho_1)\theta_{2}^\mu+\rho_2 \theta_{1}^\mu\right)+\beta_{1}^C(\rho_2+\alpha_{2}^C(\rho_1-1))}{\spio~\spit}\nonumber\\&&\hspace*{-2.5cm} -\frac{E[\mathcal{H}_{1k}^\mu|{\bf X}^\mu_{k}={\bf x}^\mu]}{\mu ~\spio~\spit} \nonumber\ +\frac{(\rho_2+\alpha_{2}^C\rho_1)s_{1}^\mu+\beta_{2}^C}{\spit} \\
 && \hspace*{-2.5cm} -\frac{E[\mathcal{H}_{2k}^\mu|{\bf X}^\mu_{k}={\bf x}^\mu]}{\mu~\spit}, \nonumber  \\
&& \hspace{-30mm} \mbox{where, } \spit= (1-\rho_2)-\alpha_{2}^C\rho_1. \label{stopping4}
\end{eqnarray}
\end{small}%
Observe that $\spit > 0$. 
Now substituting  the above value of $E[M^\mu_{2k}|{\bf X}^\mu_{k}={\bf x}^\mu]$ in \eqref{Froster1}, we have

\vspace{-2mm}
\begin{small}
        \begin{eqnarray}
    \Delta V({\bf x}^\mu)\hspace{-2mm}&\hspace{-2mm} =\hspace{-2mm}& \hspace{-2mm}\rho_2(s_1^\mu+ s_2^\mu)+\frac{(\rho_1+\rho_2-1)((\rho_2+\alpha_{2}^C\rho_1)s_{1}^\mu+\beta_{2}^C)}{\spit}\nonumber\\&& \hspace*{-0.95cm}+(\rho_1+\rho_2-1)\frac{\splitdfrac{(1-\alpha_{1}^C\alpha_{2}^C)\left((1-\rho_1)\theta_{2}^\mu+\rho_2 \theta_{1}^\mu\right)}{+\beta_{1}^C(\rho_2+\alpha_{2}^C(\rho_1-1))}}{\spio~\spit}\nonumber\\&&\hspace*{-0.95cm} +\frac{(1-\rho_1-\rho_2) E \left [ \left . \mathcal{H}_{1k}^\mu + \spio \mathcal{H}_{1k}^\mu \right  |{\bf X}^\mu_{k}={\bf x}^\mu \right ]}{\mu ~\spio~\spit} \nonumber\\
    &&\hspace*{-1.15cm}\stackrel{L}{=}(\rho_1+\rho_2-1)\frac{\splitdfrac{(1-\alpha_{1}^C\alpha_{2}^C)((1-\rho_1)(\theta_{2}^\mu-\theta_{2}^*)}{+\rho_2 (\theta_{1}^\mu-\theta_{1}^*))}}{\spio~\spit}\nonumber\\&&\hspace*{-0.95cm} +\frac{(1-\rho_1-\rho_2) E \left [ \left . \mathcal{H}_{1k}^\mu + \spio \mathcal{H}_{1k}^\mu \right  |{\bf X}^\mu_{k}={\bf x}^\mu \right ]}{\mu ~\spio~\spit} ;
    \label{Stability31}
\end{eqnarray}
\end{small}%

by Lemma \ref{Lemma_equlity_L},  the equality `$L$' is true using a unique $(\theta_1^*, \theta_2^*)$ that satisfy \eqref{stabletheta1} and \eqref{stabletheta2}. 
{ By assumption {\bf R}.1 the magnitude of the last two terms in the \eqref{Stability31} can be bounded by $\Omega $
using some $\Omega < \infty$, Now choose $\Upsilon_i,~i=1,2$ big enough such that 

\vspace*{-0.5cm}
$$
 \small\Omega  <    \frac{(1-\alpha_{1}^C\alpha_{2}^C)(1-\rho_1-\rho_2)  ((1-\rho_1)\Upsilon_2+\rho_2 \Upsilon_1)}{\spio~\spit} 
$$
possible, as
 $(1-\alpha_{1}^C\alpha_{2}^C)>0$ and $\rho_1+\rho_2 <1.$
Now define 
\begin{eqnarray}
     \mathscr{C}=\{{\bf x}^\mu:\theta_i^\mu-\theta_i^* \le \Upsilon_i,~i=1,2  \}, \label{Eqn_small_set}
\end{eqnarray}%
which has finitely many values   as the  state space 
 $$
  \left \{ (\theta_1^\mu, \theta_2^\mu): \theta_i^\mu  = \frac{1}{\mu} z_i, z_i \in \{0, 1, 2, \cdots\} \mbox{ for each } i \right  \} $$%
 is countable.
 Hence, for any $\mu$ and $\forall {\bf x}^\mu\notin \mathscr{C}$, we can have from \eqref{Stability31} 
    \begin{eqnarray*}
    \Delta V({\bf x}^\mu)&=&\frac{\splitdfrac{(1-\alpha_{1}^C\alpha_{2}^C)(\rho_1+\rho_2-1)}{((1-\rho_1)(\theta_2^\mu-\theta_2^*)+\rho_2 (\theta_1^\mu-\theta_1^*))}}{\spio~\spit}\nonumber\\\\
    &&  +\frac{\splitdfrac{(1-\rho_1-\rho_2)} {E \left [ \left . \mathcal{H}_{1k}^\mu + \spio \mathcal{H}_{1k}^\mu \right  |{\bf X}^\mu_{k}={\bf x}^\mu \right ]}}{\mu ~\spio~\spit} 
    \\ & < &  -\zeta
\end{eqnarray*} 
 
if $\alpha_{1}^C\alpha_{2}^C < 1$ and $\rho_1+\rho_2<1$. Here $\zeta>0$ is a real constant. Therefore, by Foster stability criteria (see Thorem 2.6.4 in \cite{Lyapunov}) we can say that the Markov chain in random system $\{{\bf X}_k^\mu\}_{k\in \mathbb{N}}$ is stable if  $\alpha_{1}^C\alpha_{2}^C<1$} and $\rho_1+\rho_2<1$.

\eop

\begin{lem}
\label{Lemma_equlity_L}
    The equality `$L$' in \eqref{Stability31} is satisfied. 
 \end{lem}
 {\proof} To begin with, we first    obtain 
the solution  $(\theta_1^*, \theta_2^*)$ of the following pair of equations
where (coefficients $a_{mn}$ for different $m,n$  are in \eqref{stabletheta2} at the top of the this page): 
\begin{eqnarray}
a_{11}\theta^*_{1}+a_{12}\theta^*_{2} =a_{13} \mbox{ and, }
a_{21}\theta^*_{1}+a_{22}\theta^*_{2} =a_{23}; \label{stabletheta1}
\end{eqnarray}
The solution exists as the determinant 
\begin{eqnarray}
    d=\frac{(1-\rho_1-\rho_2)(1 - \alpha_{1}^C\alpha_{2}^C)}{\spio\spit} \ne 0, \label{stabilitydit}
\end{eqnarray}
 because   $\alpha_{1}^C\alpha_{2}^C\neq 1$ and $\rho_1+\rho_2 \neq 1.$
Now  combine $(\theta^*_{1},\theta^*_{2})$  as below and observe that equality $L$ is satisfied:
\begin{eqnarray}
    &&\frac{(1-\rho_1-\rho_2)(1-\alpha_{1}^C\alpha_{2}^C)\left((1-\rho_1)\theta_{2}^*+\rho_2 \theta_{1}^*\right)}{\spio~\spit}\nonumber\\&&=\rho_2(s_1^\mu+ s_2^\mu)+\frac{(\rho_1+\rho_2-1)((\rho_2+\alpha_{2}^C\rho_1)s_{1}^\mu+\beta_{2}^C)}{\spit}\nonumber\\&&+\frac{(\rho_1+\rho_2-1)\beta_{1}^C(\rho_2+\alpha_{2}^C(\rho_1-1))}{\spio~\spit}\label{stability2}.  
\end{eqnarray} 
\eop

{\bf Proof of Theorem 
\ref{thm_palmmoments}:}\label{Proof_thm_2} 
In the following we compute the  expectations  $E[\Theta_{ik}^\mu]$ and   $E[\Psi_{ik}^\mu]$. 
Taking  expectations of 
  equation \eqref{visiting1}    (after normalizing with $\mu$) and rearranging we have (also using    \eqref{Eqn_Lam_Gam}, \eqref{Eqn_Dik}, \eqref{Queue1} and \eqref{Queue2}): 
    \begin{eqnarray}
E[M^\mu_{1k}] =\frac{E[\Theta_{1k}^\mu]+\alpha_{1}^C E[\Theta_{2k}^\mu]+\beta_{1}^C}{\spio}-\frac{E[\mathcal{H}_{1k}^\mu]}{\mu\spio},\label{stopping5}
 \end{eqnarray}
 
These steps are as in the proof of Theorem \ref{thm_palmstability}. 
Now using  equations \eqref{visiting2}  and \eqref{stopping5} and proceeding as above,  
\begin{eqnarray}
 E[M^\mu_{2k}] &=& \frac{(1-\alpha_{1}^C\alpha_{2}^C)\left((1-\rho_1)E[\Theta_{2k}^\mu]+\rho_2 E[\Theta_{1k}^\mu]\right)}{\spio~\spit}\nonumber\\&& +\frac{\mu\beta_{1}^C(\rho_2+\alpha_{2}^C(\rho_1-1))-E[\mathcal{H}_{1k}^\mu]}{\mu ~\spio~\spit} \nonumber\\&&+\frac{(\rho_2+\alpha_{2}^C\rho_1)s_{1}^\mu+\beta_{2}^C}{\spit}-\frac{E[\mathcal{H}_{2k}^\mu]}{\mu~\spit}.\label{stopping6}
\end{eqnarray}

By normalizing the equations \eqref{THETAi} for $i=1,2$ with $\mu$
   and  taking the  expectation, we obtain the following using \eqref{Eqn_Lam_Gam},\vspace*{-0.3cm}
   
   \begin{small}
       \begin{eqnarray*}
   E\left[\Theta_{i(k+1)}^\mu\right] &=& E[\Theta_{ik}^\mu]+\rho_i\left(E[M^\mu_{1k} + [M^\mu_{2k}]+s_1^\mu+ s_2^\mu\right)\\ && -E[M^\mu_{ik}].\label{Implicit1}
 \end{eqnarray*}
    \end{small}
    
Using the   values of $E[M^\mu_{ik}]$ in the above equation, we have \vspace*{-0.3cm}

   \begin{small}
    \begin{eqnarray*}
    E[\Theta^\mu_{1(k+1)}] &=&
    (1+a_{11}) E[\Theta^\mu_{1k}]+a_{12}E[\Theta^\mu_{2k}]-a_{13} + \frac{r_{1k}}{\mu}\\[2.5pt]
    E[ \Theta^\mu_{2(k+1)}] &=&a_{21} E[\Theta^\mu_{1k}]+(1+a_{22})E[\Theta^\mu_{2k}]-a_{23}+\frac{r_{2k}}{\mu}
     \label{moment1}
\end{eqnarray*}
\end{small}%
where the terms $\nicefrac{r_{ik}}{\mu}$ are constructed using $E[\mathcal{H}_{ik}^\mu]$ of \eqref{stopping5}-\eqref{stopping6}. 

Define $m_{ik} := E[\Theta_{ik}^\mu - \theta_i^*]$ for each $i, k$ and using \eqref{stabletheta1} in above equations we have matrix equation:

\vspace*{-0.3cm}
\begin{small}
    \begin{eqnarray*}
    m_{1(k+1)} &=&
    (1+a_{11}) m_{1k}+a_{12}m_{2k} + \frac{r_{1k}}{\mu}\\[2.5pt]
    m_{2(k+1)} &=&a_{21} m_{1k}+(1+a_{22})m_{2k}+\frac{r_{2k}}{\mu}
     \label{moment2}
\end{eqnarray*}
\end{small}

Define $\overline{\bf m}_{k}:=[m_{1k},m_{2k}]^{\mathsf {tr}}$,  $\overline{\bf r}_{k}:=[r_{1k},r_{2k}]^{\mathsf{tr}}$ for each $k$. Therefore, 
\begin{equation*}
    \overline{\bf m}_{k+1}= {\bf A}\overline{\bf m}_{k}+\frac{1}{\mu}\overline{\bf r}_{k}, \quad {\bf A}= \left[ \begin {array}{cc} 
1+a_{11}& a_{12}   \\
a_{21} & 1+a_{22}
\end {array} \right].
\end{equation*}

Hence, we have 
$$\overline{\bf m}_{k+1}= {\bf A}^k\overline{\bf m}_{1}+\frac{1}{\mu}\sum\limits_{l=1}^k {\bf A}^{k-l}\overline{\bf r}_{l}.$$
Now, as ${\bf Id}-{\bf A}$ has the determinant $d$ given in \eqref{stabilitydit}, which is non-zero from the condition $\alpha_{1}^C\alpha_{2}^C\neq 1$ and $\rho_1+\rho_2 \neq 1.$  Thus ${\bf A}^k \to {\bf O}$, the zero matrix.  

Also, from ${\bf R.}1$, we can have $|r_{ik}|<r_i$ uniformly in $k$, for each $i$, where $r_i$ is some constant. Therefore, $\overline{\bf r}_{k}$ is bounded uniformly in $k$. Thus $E[\Theta_{ik}^\mu] \to \theta_i^* + \nicefrac{r_i}{\mu} \mbox{ for each }  i \in {1, 2}$, when $k \to \infty$. Substituting this in \eqref{stopping5}-\eqref{stopping6} we obtain the rest of the proof as, $E[\Psi^\mu_k] = E[M_{1k}^\mu+M_{2k}^\mu]+s_1^\mu+s_2^\mu.$

\eop

\begin{table}[ht]
\scalebox{0.99}{
\begin{tabular}{|c|c|}
\hline
{\bf Notation}             & {\bf Description}                                                                                          \vspace{1mm}\\ \hline
$\lambda, \mu, Q_i, \rho_i$             & Arrival rate, Service rate,  i-th Queue, traffic load at $Q_i$                                                \vspace{0.08mm}\\ \hline

  $\Pi_{iB},\Pi_{iC}$                  &  $\substack{\mbox{Switching functions relevant while visiting $Q_i$ } \\ \mbox{
  }}$                         \\  \hline \vspace{0.08mm}
 $\substack{\alpha_{i}^B,\beta_i^B, \\\alpha_{i}^{C}, \beta_{i}^C}$                  &  $\substack{\mbox{Switching parameters relevant while visiting $Q_i$ } \\ 
 \mbox{  }
 }$                         \\ \hline

$S_i^\mu$                  &  Switching time from $Q_i$ to $Q_{-i}$ in $\mu$-system                                               \vspace{0.08mm} \\ \hline
$\Lambda_i^\mu$                  &  Arrival process for $Q_i$ in $\mu$-system                                                \vspace{0.08mm}\\ \hline 
$\Gamma_{ik}^\mu$                  &  $\substack{\mbox{Uninterrupted service process at} \\ k \mbox{-th visit for $Q_i$ in $\mu$-system}}$ \vspace{0.1mm}\\ \hline  

$B^\mu_{ik}, C^\mu_{ik}$             & $\substack{\mbox{Length of Beginning phase, Concluding phase}\\ \mbox{for $k$-th visit of $Q_i$ in $\mu$-system}  }$                            \vspace{1mm}\\ \hline
$M^\mu_{ik}$             & {Total visit time for $k$-th visit of $Q_i$     in $\mu$-system}                          \vspace{0.08mm}\\ \hline
$N_i^\mu(t)$                & $\substack{\mbox{Number of waiting customers at $Q_i$}\\ \mbox{at time $t$ in $\mu$-system}}$ \vspace{0.1mm}\\ \hline
$\Psi_k^\mu$                &                                        Duration of $k$-th cycle in $\mu$-system    \vspace{0.08mm}\\ \hline
$\Phi_k^\mu$                &   Starting epoch of $k$-th cycle  in $\mu$-system                                      \vspace{0.08mm}\\ \hline
$N_{ik}^\mu$                &   Number of waiting customers for  $Q_i$ at $\Phi^\mu_k$                                     \vspace{0.08mm}\\ \hline
$\Theta_{ik}^\mu$                &  $\substack{\mbox{Normalized  number of waiting }\\ \mbox{customers for  $Q_i$ at $\Phi^\mu_k$}}$                                      \vspace{0.1mm}\\ \hline
${\bf X}_{k}^\mu$                &  $\substack{\mbox{Joint Palm Markov chain  (discrete-time) for normalized}\\\mbox{  number of customers at both queues in $\mu$-system}} $                              \vspace{0.1mm}\\ \hline
${\bf X}^\mu(t)$                & $\substack{\mbox{Joint stochastic process for normalized number of}\\\mbox{ customers at both queues in $\mu$-system}} $                           \vspace{0.1mm}\\ \hline

$R_m^\mu$  & Regeneration epochs for ${\bf X}^\mu(t)$                          \vspace{0.08mm}\\ \hline

${\cal T}_s$ & Class of stable switching policies  \vspace{0.08mm}\\ \hline
${\cal P}_s$ & $\substack{\mbox{Subclass of stable switching policies  }\\\mbox{ generating Pareto frontier}}$\vspace{0.09mm}\\ \hline
\end{tabular}}
\caption{Summary of Notations}
\end{table}

\end{document}